\documentclass[11pt]{article}
\newtheorem{teor}{Theorem}[section]
\newtheorem{lema}[teor]{Lemma}
\newtheorem{prop}[teor]{Proposition}
\newtheorem{corol}[teor]{Corollary}
\newtheorem{defin}[teor]{Definition}
\newtheorem{rem}[teor]{Remark}

\usepackage{amssymb,amsmath}
\usepackage{hyperref}
\usepackage[english]{babel}

\title{On the geometry of almost Golden Riemannian manifolds}
\author{Fernando Etayo\footnote{Departamento de Matem\'{a}ticas, Estad\'{\i}stica y Computaci\'{o}n. Facultad de Ciencias.  Universidad de Cantabria. Avda. de los Castros, s/n, 39071 Santander, SPAIN. e-mail: etayof@unican.es},\, Rafael Santamar\'{\i}a\footnote{Departamento de Matem\'{a}ticas. Escuela de Ingenier\'{\i}as Industrial e Inform\'{a}tica. Universidad de Le\'{o}n. Campus de Vegazana, 24071 Le\'{o}n, SPAIN. e-mail: rsans@unileon.es}\, and Abhitosh Upadhyay\footnote{Harish-Chandra Research Institute, HBNI, Chhatnag Road, Jhunsi, Allahabad, 211019, INDIA. e-mail: abhi.basti.ipu@gmail.com, abhitoshupadhyay@hri.res.in}}

\date{}
\begin{document}
\maketitle
\begin{abstract}
An almost Golden Riemannian structure $(\varphi ,g)$ on a manifold is given by a tensor field $\varphi $ of type (1,1) satisfying the Golden section relation $\varphi ^{2}=\varphi +1$, and a pure Riemannian metric $g$, i.e., a metric satisfying $g(\varphi  X,Y)=g(X,\varphi  Y)$. We study connections adapted to such a structure, finding two of them, the first canonical and the well adapted, which measure the integrability of $\varphi $ and the integrability of the $G$-structure corresponding to $(\varphi ,g)$.

\end{abstract}



\vspace{7mm}

\section{Introduction} The Golden section or Golden mean $\phi$ is the positive root of the polynomial equation $x^2-x-1=0$; i.e., $\phi=\frac{1+\sqrt{5}}{2}$. The negative root of the previous equation, usually denoted by $\bar\phi$, satisfies $\bar\phi=\frac{1-\sqrt{5}}{2}=1-\phi$. In the last years the Golden mean can be found in many areas of mathematical and physical research (see, {\em e.g.}, \cite{elnaschie} and \cite{stakhov} and the references therein). As generalization of the Golden mean appear the metallic means (see \cite{spinadel}), which are the positive root of the equation $x^2-px-q=0$, where $p,q$ are positive integers.

Analogously to almost product and almost complex structures on a differentiable manifold, in \cite{crasmareanu}, Crasmareanu and Hre\c{t}canu introduced and studied the Golden structures. All of them are examples of polynomial structures of degree 2. Goldberg and Yano defined the polynomial structures of degree $d$ on a manifold $M$ as tensor fields $J$ of type $(1,1)$ with constant rank such that $d$ is the smallest integer which $J^d, J^{d-1}, \ldots, J, Id$ are not independent, where $Id$ is the identity tensor field (see \cite{goldberg}). Then there exist $a_{d-1}, \ldots, a_1, a_0 \in \mathbb{R}$ such that $J^d+a_{d-1} J^{d-1}+\ldots + a_1 J +a_0 Id=0$. An almost product (resp. almost complex) structure $J$ satisfies $J^2-Id=0$ (resp. $J^2+Id=0$). In this case $(M,J)$ is called almost product (resp. almost complex) manifold. Now we precise the notion of almost Golden structure.

\begin{defin}[{\cite[Defin. 1.2]{crasmareanu}}] Let $M$ be a manifold. A polynomial structure $\varphi$ of degree 2 on $M$ satis\-fying  $\varphi^2=\varphi+Id$
is called an almost golden structure. In this case $(M,\varphi)$ is an almost Golden manifold.
\end{defin}

Some authors later study metallic structures which are polynomials structures satisfying $\varphi^2=p\varphi+qId$, where $p,q$ are positive integers.
Then, the first example of these structures are the Golden ones.

Unlike the paper \cite{crasmareanu} we use the word ``almost'' to general Golden structures on a manifold. We say that a manifold $M$ is a Golden manifold if it has an integrable almost Golden structure.

Recall that a polynomial structure $J$ is integrable if  the Nijenhuis tensor $N_{J}$ vanishes, where
\begin{equation}
\label{eq:nijenhiis}
N_{J}(X,Y)=J^2[X,Y]+[JX, JY]-J [JX, Y]-J[X, JY], \quad \forall X, Y \in {\mathfrak X} (M).
\end{equation}

Many works about polynomial structures are devoted to study compatible Riemannian metrics and adapted or natural covariant derivatives to a polynomial structure. Let $J$ be a polynomial structure on $M$. A  Riemannian metric $g$ is compatible with $J$ if satisfies
\begin{equation}
\label{eq:g-compatible}
g(JX,Y)=g(X,JY), \quad \forall X, Y \in {\mathfrak X} (M).
\end{equation}
The most known manifolds in these conditions are almost Riemannian product with null trace manifolds. This kind of manifolds are a particular class of $(J^2=\pm1)$-metric ones (see \cite{brno}). If the metric $g$ satisfies the condition (\ref{eq:g-compatible}) then $g$ is called a pure metric respect to the polynomial structure $J$ (see, {\em e.g.}, \cite{gezer} and \cite{salimov}).

One says that a covariant derivative $\nabla$ on $M$ is adapted to $J$ if $\nabla J=0$. As is well known, if $\nabla$ is a torsion-free covariant law then one has
\begin{equation}
N_J(X,Y) =(\nabla_X J) JY + (\nabla_{JX} J) Y - (\nabla_Y J) JX -(\nabla_{JY} J) X, \quad \forall X, Y \in {\mathfrak X} (M).
\label{eq:nijenhuis-nabla}
\end{equation}
Then the existence of symmetric adapted derivation laws to a polynomial structure is a sufficient condition to assure its integrability.
\medskip

In \cite{crasmareanu}, the authors also studied almost Golden structures which admit compatible Riemannian metrics.

\begin{defin}[{\cite[Defin. 5.1]{crasmareanu}}]
\label{teor:agr-condition}
 Let $(M,\varphi)$ be an almost Golden structure and let $g$ be a Riemannian metric on $M$. We say that $(\varphi,g)$ is an almost Golden Riemannian structure on $M$ if
\begin{equation}
\label{eq:agr-condition}
g(\varphi X, Y)= g(X, \varphi Y), \quad \forall X, Y \in {\mathfrak X} (M).
\end{equation}
In this case $(M,\varphi,g)$ is called an almost Golden Riemannian manifold.
\end{defin}
Note that the condition \eqref{eq:agr-condition} is equivalent to the following one
\[
g(\varphi X, \varphi Y) =g(\varphi X, Y)+ g(X, Y), \quad \forall X, Y \in {\mathfrak X} (M).
\] 
However in the case of almost product Riemannian manifolds the condition (\ref{eq:g-compatible}) is equivalent to
\begin{equation}
\label{eq:g-almostproduct}
g(JX, JY)= g(X, Y), \quad \forall X, Y \in {\mathfrak X} (M).
\end{equation}

The main technique to study the geometry of an almost Golden manifold in \cite{crasmareanu} is based on the use of the corresponding almost product structure. The authors show that almost Golden and almost product structures are closely related by the following

\begin{prop}[{\cite[Theor. 1.1]{crasmareanu}}]
\label{teor:equivalence-structures}
Let $M$ be a manifold.  
\begin{enumerate}

\item[$i)$] Let $\varphi$ be an almost Golden structure on $M$. Then 
\begin{equation}
\label{eq:J-induced}
J_{\varphi} = \frac{1}{\sqrt{5}}(2\varphi-Id)
\end{equation}
is an almost product structure on $M$. We say that $J_{\varphi}$ is the almost product structure induced by $\varphi$.

\item[$ii)$] Let $J$ be an almost product structure on $M$. Then  $\varphi_{J} = \frac{1}{2}(Id+\sqrt{5}J)$ is an almost Golden structure on $M$. We say that $\varphi_{J} $ is the almost Golden induced by $J$.
\end{enumerate}
Thus there exist an $1\colon1$ correspondence between almost Golden structures and almost product structures on $M$ because $\varphi_{J_{\varphi}}=\varphi$ and $J_{\varphi_{J}}=J$.
\end{prop}

As direct consequence of the introduction of the almost product structure induced by the almost Golden structure one has the following  

\begin{lema}
\label{teor:nijenhuis-varphi-induced}
Let $(M,\varphi)$ be an almost Golden manifold. Then
\[
N_{J_{\varphi}}(X,Y)=\frac{4}{5} N_{\varphi}(X,Y), \quad \forall X, Y \in {\mathfrak X} (M).
\]
\end{lema}

{\bf Proof.} By straightforward calculations from the definitions of the almost product structure $J_{\varphi}$ induced by $\varphi$ and the Nijenhuis tensors of $\varphi$ and $J_{\varphi}$ (see (\ref{eq:nijenhiis}) and (\ref{eq:J-induced})). $\blacksquare$
\medskip 

Then it is obvious the equivalence between the integrability of the structures $\varphi$ and $J_{\varphi}$.
\medskip

Starting from the results of Crasmareanu and Hre\c{t}canu we will continue the study of this kind of manifolds using the closely relation between almost Golden and almost product structures. We want to go further into the relation specified in Proposition \ref{teor:equivalence-structures}. Then it is necessary to link the $G$-structures defined over $M$ by the two polynomial structures. 

Bearing in mind the Definition \ref{teor:agr-condition}, Proposition \ref{teor:equivalence-structures} and condition (\ref{eq:g-compatible}) it easy to prove the following equivalencies 
\begin{eqnarray*}
g(\varphi X, Y) = g(X, \varphi Y) &\Longleftrightarrow& g(J_{\varphi} X, Y) = g(X, J_{\varphi} Y), \\
g(JX, Y) = g(X, JY) &\Longleftrightarrow& g(\varphi_J X, Y) = g(X, \varphi_J Y), \quad \forall X, Y \in {\mathfrak X} (M),
\end{eqnarray*}
i.e., the Riemannian metric $g$ is pure with respect to the almost Golden structure $\varphi$ (resp. $\varphi_J$)  if and only  $g$ is pure respect to the almost product structure $J_{\varphi}$ (resp. $J$) too. Then the next step will be study the $G$-structures defined by both structures in the case when they admit a compatible Riemannian metric. By means of these $G$-structures we will focus on almost Golden Riemannian manifolds. We will show the set of natural derivation laws of this kind of manifolds, which are the covariant derivatives that simultaneously parallelize the Riemannian metric and the almost Golden structure (see Definition \ref{eq:natural-conection-goldenR}). Previously we will address the set of natural connections of an almost Golden structure 
highlighting its relation with the set of natural connections of the almost product structure induced by the almost Golden structure.

We will prove that there exist two distinguished natural derivation laws on an almost Golden Riemannian manifold, the first canonical and the well adapted connections, which characterize  the integrability of the almost Golden structure (Theorem \ref{teor:toro-varphi}) and the $G$-structure defined over an almost Golden Riemannian manifold (Proposition \ref{teor:g-structure-integrable}). We will conclude our study analyzing the coincidence of both linear connections and the Levi Civita connection  of the Riemannian metric. The main difference between Schouten and Vr\u{a}nceanu connections in the case of almost Golden manifolds and the previous ones in the case of almost Golden Riemannian manifolds is that the first canonical and the well adapted connections can be uniquely defined.

We will consider smooth manifolds and operators being of class $C^{\infty}$. As in this introduction, $\mathfrak{X}(M)$ denotes the module of vector fields of a manifold $M$. We denote by ${\mathcal T}^p_q(M)$ the module of tensor fields of type $(p,q)$ of a manifold $M$. The general (resp. orthogonal) linear group will be denoted as usual by $GL(n,\mathbb{R})$ (resp. $O(n,\mathbb{R})$). 

\section{$G$-structures defined over almost Golden Riemannian manifolds}

Now we analyze the relation between the $G$-structures defined over a manifold by almost Golden and almost product structures. Let $(M,\varphi)$ be an almost Golden structure. For every $p\in M$ the tangent vector space $T_p(M)$ splits as follows 
\[
T_p(M)=(D_{\phi})_p\oplus (D_{\bar\phi})_p = T_p^+(M)\oplus T^-_p(M),
\]
where
\begin{eqnarray*}
(D_{\phi})_p=\{X \in T_p(M)\colon \varphi_pX=\phi X\}, &&(D_{\bar\phi})_p=\{X \in T_p(M)\colon \varphi_pX=\bar\phi X\}, \\
T_p^+(M)=\{X \in T_p(M)\colon (J_\varphi)_p  X=X\}, && T_p^-(M)=\{X \in T_p(M)\colon (J_\varphi)_p  X=-X\}.
\end{eqnarray*}
The relation between these vector subspaces of $T_p(M)$ is the following
\begin{eqnarray*}
&\varphi_p X = \phi X \Longleftrightarrow \varphi_p X = \frac{1+\sqrt{5}}{2} X \Longleftrightarrow \frac{1}{\sqrt{5}} (2 \varphi_p X-X) =X \Longleftrightarrow (J_\varphi)_p X =X,& \\
&\varphi_p X = \bar\phi X \Longleftrightarrow \varphi_p X = \frac{1-\sqrt{5}}{2} X \Longleftrightarrow \frac{1}{\sqrt{5}} (2 \varphi_p X-X) =-X \Longleftrightarrow (J_\varphi)_p X =-X, & \forall p \in M, \forall X \in T_p(M).
\end{eqnarray*}
If one denotes $T^+(M)=\ker (J_{\varphi}-Id)$ and $T^-(M)=\ker (J_{\varphi}+Id)$, then 
the above equivalences show that
\[
(D_{\phi})_p=T_p^+(M), \quad  (D_{\bar\phi})_p=T_p^-(M), \quad \forall p \in M,
\]
and then the corresponding distributions coincide: $D_{\phi}=T^+(M)$, $D_{\bar\phi}=T^-(M)$. Then one can claim

\begin{prop}
\label{teor:gstructure-golden}
Let $(M,\varphi)$ be an $n$-dimensional almost Golden manifold. The polynomial structure $\varphi$ defines the following $G_{\varphi}$-structure  over $M$
\[
\mathcal G_{\varphi} =\bigcup_{p\in M} 
\left\{ (X_1, \ldots, X_r, Y_1, \ldots ,Y_s) \in \mathcal F_p (M) \colon 
X_1, \ldots , X_r \in (D_{\phi})_p, Y_1, \ldots  , Y_s \in (D_{\bar\phi})_p
\right\},
\]
whose structural group is $G_{\varphi}=GL(r;\mathbb{R})\times GL(s;\mathbb{R})$, where $\mathcal F(M)$ is the principal bundle of linear frames of $M$ and $r=\dim D_{\phi}$, $s=\dim D_{\bar\phi}$, $r+s=n$.
\end{prop}

{\bf Proof.} Trivial. $\blacksquare$
\medskip

Therefore the $G_{\varphi}$-structure defined by $\varphi$ over $M$ is the same $G$-structure defined by $J_\varphi$.
\medskip

The above identities justify the existence of the $1:1$ correspondence between almost Golden and almost product structures on a manifold $M$ showed in Proposition \ref{teor:equivalence-structures}.
\medskip

The next step in our study is the $G$-structure defined by an almost Golden Riemannian structure. Proposition \ref{teor:gstructure-golden} carries naturally to the following one.

\begin{prop}
Let $(M,\varphi,g)$ be an $n$-dimensional almost Golden Riemannian manifold. The $G_{(\varphi,g)}$-struc\-tu\-re over $M$ defined by $(J,g)$ is
\[
\mathcal G_{(\varphi,g)} =\bigcup_{p\in M} 
\left\{ (X_1, \ldots, X_r, Y_1, \ldots ,Y_s) \in \mathcal F_p (M) \colon 
\begin{array}{c}
X_1, \ldots , X_r \in (D_{\phi})_p, Y_1, \ldots  , Y_s \in (D_{\bar\phi})_p,\\
g(X_i, X_j)=\delta_{ij}, i, j=1, \ldots,  r, \\
g(Y_i, Y_j)=\delta_{ij}, i, j=1, \ldots, s, \\
g(X_i, Y_j)=0, i,=1, \ldots r, j=1,\ldots, s \\
\end{array}
\right\},
\]
whose structural group is $G_{(\varphi,g)}=O(r;\mathbb{R})\times O(s;\mathbb{R})$
\end{prop}

{\bf Proof.} Trivial. $\blacksquare$
\medskip

Then it is obvious that the $G_{(\varphi,g)}$-structure coincides  with the $G$-structure defined by the almost product Riemannian structure $(J_{\varphi},g)$ over $M$.
\medskip

The previous results consolidate the approach to the study of almost Golden manifolds initiated in \cite{crasmareanu} by using the almost product structure induced by the almost Golden one. In this paper, we explore the geometry of an almost Golden Riemannian manifold by means of  the induced almost product Riemann structure.

\section{Adapted connections to almost Golden structures}

The set of all linear connections whose derivation law is adapted to $\varphi$ can be found in \cite[Theor 5.1]{crasmareanu}. We will analyze later this result more deeply. First, we relate derivation laws adapted to an almost Golden structure with derivation laws adapted to an almost product structure.

\begin{lema}
\label{teor:equivalence-adapted-connections}
Let $(M,\varphi)$ be an almost Golden manifold and let $\nabla$ be a derivation law on $M$. Then $\nabla$ is adapted to $\varphi$ if and only if $\nabla$ is adapted to the almost product structure $J_{\varphi}$ induced by $\varphi$.

\end{lema}

{\bf Proof.} Recall that if $J$ is a tensor field of type $(1,1)$ then $\nabla J=0$ means $J \nabla_X Y = \nabla_X JY$ for every $X$, $Y$, vector fields on $M$.

Suppose that $\nabla\varphi=0$ then
\[
\nabla_X J_{\varphi} Y = \frac{2}{\sqrt 5} \nabla_X \varphi Y - \frac{1}{\sqrt 5}  \nabla_X Y = \frac{2}{\sqrt 5} \varphi(\nabla_X Y)- \frac{1}{\sqrt 5} \nabla_X Y = J_{\varphi} (\nabla_X Y), \quad \forall X, Y \in {\mathfrak X} (M), 
\]
thus $\nabla J_{\varphi}=0$.  Analogously if $\nabla J_{\varphi}=0$ then 
\[
\nabla_X \varphi Y = \frac{1}{2} \nabla_X Y + \frac{\sqrt{5}}{2} \nabla_X J_{\varphi} Y = \frac{1}{2} \nabla_X Y +    \frac{\sqrt{5}}{2} J_{\varphi} \nabla_X Y = \varphi(\nabla_X Y), \quad \forall X, Y \in {\mathfrak X} (M),
\]
therefore $\nabla \varphi=0$. $\blacksquare$
\medskip

The previous lemma shows that the shortest way to attach an adapted derivation law to an almost Golden structure is to build an adapted one to the induced almost product structure. If we denote by $J_{\varphi}^+$ and $J_{\varphi}^-$ the projectors of $TM$ over $T^+M$ and $T^-M$ respectively, starting from a derivation law $\nabla$ one can construct the following two adapted derivation laws to the almost product structure $J_{\varphi}$
\begin{eqnarray*}
\nabla^{\mathrm{sc}}_X Y&=& J_{\varphi}^+ \nabla_X J_{\varphi}^+Y + J_{\varphi}^-\nabla_X J_{\varphi}^-Y,\\
\nabla^{\mathrm{v}}_X Y&=& J_{\varphi}^+ \nabla_{J_{\varphi}^+X} J_{\varphi}^+Y + J_{\varphi}^-\nabla_{J_{\varphi}^-X} J_{\varphi}^-Y +J_{\varphi}^+[J_{\varphi}^-X, J_{\varphi}^+Y]+J_{\varphi}^-[J_{\varphi}^+X,J_{\varphi}^-Y], \quad \forall X, Y \in {\mathfrak X} (M),
\end{eqnarray*}
called the Schouten and Vr\u{a}nceanu type connections (see {\em e.g.}; \cite{ianus} and the references therein).

An almost product structure is a particular case of an $\alpha$-structure, which is a way of generalizing simultaneously almost product and almost complex structures on a manifold (see \cite[Sec. 2]{brno}). In that paper, the authors study linear connections whose covariant derivatives are adapted to an $\alpha$-structure. By this reason they introduce the so-called $\nabla^{0}$ type connections starting from a derivation law $\nabla$ as in the previous cases
\[
\nabla^{0}_X Y = \nabla_X Y -\frac{1}{2}( \nabla_X J_{\varphi}) J_{\varphi}Y, \quad \forall X, Y \in {\mathfrak X} (M),
\]
which is also a natural derivation law of $J_{\varphi}$. It is easy to prove that $\nabla^{0}$ and Schouten type connections are the same in the case of the polynomial structure $J_{\varphi}$. Indeed, taking into account that
\[
Id=J_{\varphi}^+ + J_{\varphi}^-, \quad J_{\varphi}=J_{\varphi}^+ - J_{\varphi}^-,
\]
then for every $X, Y$ vector fields on $M$ one has
\begin{eqnarray*}
\nabla^{0}_X Y &=& \nabla_X Y -\frac{1}{2}( \nabla_X J_{\varphi}) J_{\varphi}Y= \nabla_X Y  -\frac{1}{2} \nabla_X Y   +\frac{1}{2} J_{\varphi} \nabla_X J_{\varphi}Y =\frac{1}{2} \nabla_X Y   +\frac{1}{2} J_{\varphi} \nabla_X J_{\varphi}Y\\
&=& \frac{1}{2} \left( J_{\varphi}^+\nabla_X J_{\varphi}^+Y + J_{\varphi}^+\nabla_X J_{\varphi}^-Y + J_{\varphi}^-\nabla_X J_{\varphi}^+Y + J_{\varphi}^-\nabla_X J_{\varphi}^-Y   \right)\\
&+& \frac{1}{2} \left( J_{\varphi}^+\nabla_X J_{\varphi}^+Y - J_{\varphi}^+\nabla_X J_{\varphi}^-Y - J_{\varphi}^-\nabla_X J_{\varphi}^+Y + J_{\varphi}^-\nabla_X J_{\varphi}^-Y   \right)\\
&=&J_{\varphi}^+ \nabla_X J_{\varphi}^+Y + J_{\varphi}^-\nabla_X J_{\varphi}^-Y =\nabla^{\mathrm{sc}}_X Y.
\end{eqnarray*}

These adapted connections to the almost Golden structure  can be used to study $D_{\phi}$ and $D_{\bar\phi}$ (see \cite[Prop. 4.2]{crasmareanu} and \cite[Prop. 2.1]{ianus}). Moreover, $\nabla^{0}$ type connections, or Schouten type connections if you prefer, allow us to describe the set of natural derivation laws of the almost product structure $J_{\varphi}$ (see \cite[Prop. 2.4]{brno}).

However, the Schouten (or $\nabla^{0}$) and Vr\u{a}nceanu type connections are not uniquely defined because  they depend on an arbitrary covariant derivative.
\medskip

$\nabla^{0}$ type connections will help us to show the set of linear connections adapted to an almost Golden structure by means of the induced almost product structure. 

\begin{prop}
\label{teor:adaptedconnections-golden}
Let $(M,\varphi)$ be an almost Golden manifold and let $\nabla$ be a derivation law on $M$. The set of derivation laws adapted to $\varphi$ is 
\[
\{
\nabla^{0}+ Q \colon Q \in \mathcal{L}\},
\]
where 
\[
\mathcal{L}=\{ Q \in \mathcal T^1_2(M) \colon Q(X,J_{\varphi}Y)=J_{\varphi} Q(X,Y), \forall X, Y  \in {\mathfrak X} (M)\}.
\]
\end{prop}

{\bf Proof.} Bearing in mind that one can describe the set of natural derivation laws of $J_{\varphi}$ by means of $\nabla^{0}$ as above  (see \cite[Prop. 2.4]{brno}), Lemma \ref{teor:equivalence-adapted-connections} allows to claim that the set of adapted derivation laws to an almost Golden structure $\varphi$  is the set presented above. $\blacksquare$

\begin{rem}
Crasmareanu and Hre\c{t}canu proved that the set of adapted covariant derivatives to an almost Golden structure in \cite[Theor. 5.1]{crasmareanu} is the following
\[
\left\{\frac{1}{5}(3\nabla_X Y + 2\varphi(\nabla_X \varphi Y )-\varphi(\nabla_X Y)-\nabla_X \varphi Y)+O_{J_{\varphi}}Q(X,Y) \colon X, Y \in {\mathfrak X} (M), Q \in \mathcal T^1_2(M)\right\},
\]
where $O_{J_{\varphi}}Q(X,Y)=\frac{1}{2}(Q(X,Y)+J_{\varphi} Q(X,J_{\varphi} Y))$, for every $X$ and $Y$ vector fields on $M$, is called the associated Obata operator.
This set coincides with the previous one.  An easy calculation shows that the set $\mathcal L$ of the previous proposition is the image of $\mathcal T^1_2(M)$ by the associated Obata operator (see \cite[Prop. 2.3]{brno}). Taking into account the definition of $J_{\varphi}$ one has
\begin{eqnarray*}
\nabla^{0}_X Y &=& \nabla_X Y -\frac{1}{2}( \nabla_X J_{\varphi}) J_{\varphi}Y =  \nabla_X Y - \frac{1}{2} (\nabla_X  Y -J_{\varphi} (\nabla_X J_{\varphi} Y))\\
              &=& \frac{1}{2}\nabla_X Y +\frac{1}{10}(4\varphi(\nabla_X \varphi Y)-2\varphi(\nabla_X Y)-2\nabla_X \varphi Y+\nabla_X Y)\\
              &=&\frac{1}{5}(3\nabla_X Y + 2\varphi(\nabla_X \varphi Y )-\varphi(\nabla_X Y)-\nabla_X \varphi Y), \quad \forall X, Y \in {\mathfrak X} (M).
\end{eqnarray*}
\end{rem}

As we have pointed out in the Introduction, adapted connections are useful tools to study the integrability of a polynomial structure (see identity (\ref{eq:nijenhuis-nabla})). In the case of an almost product structure, integrability can be characterized by the existence of a symmetric adapted covariant law (see \cite[Theor. 25]{yano} and \cite[Prop. 2.13]{brno}). Then, as direct consequence of the quoted results and Lemma \ref{teor:equivalence-adapted-connections} one has

\begin{prop}
\label{teor:goldenintegrable}
Let $(M,\varphi)$ be an almost Golden manifold. The following two conditions are equivalent
\begin{enumerate}
\item[$i)$] The almost Golden structure $\varphi$ is integrable.
\item[$ii)$] There exists a torsion-free derivation law $\nabla$ such that $\nabla \varphi=0$.
\end{enumerate}
\end{prop}

\section{Adapted connections to almost Golden Riemannian structures}

We continue the study of almost Golden Riemannian manifolds starting from Proposition \ref{teor:adaptedconnections-golden}. We will study the covariant derivatives  which parallelize the almost Golden structure and the Riemannian metric.

\begin{defin} 
\label{eq:natural-conection-goldenR}
Let $(M,\varphi,g)$ be an almost Golden Riemannian manifold and let $\nabla$ be a derivation law on $M$. The derivation law $\nabla$ is said to be natural or adapted to $(\varphi,g)$ if $\nabla \varphi=0$, $\nabla g=0$.
\end{defin}

If one changes the almost Golden structure by an almost product structure in the last definition then one recovers the  well known notion of adapted derivation law on almost product Riemannian manifolds.
Thus one has the following result analogous to Lemma \ref{teor:equivalence-adapted-connections}.

\begin{prop}
\label{teor:equivalence-adapted-connections-metric}
Let $(M,\varphi,g)$ be an almost Golden Riemannian manifold and let $\nabla$ be a derivation law on $M$. The derivation law $\nabla$ is adapted to the almost Golden Riemannian structure $(\varphi,g)$ if and only if $\nabla$ adapted to the almost product Riemannian structure $(J_{\varphi},g)$.
\end{prop}

Unlike the previous case, there exists a distinguished connection on a Riemannian manifold, the Levi Civita connection, which allow to us define a distinguished natural connection on an almost Golden Riemannian manifold.

\begin{defin}Let $(M,\varphi,g)$ be an almost Golden Riemannian manifold and let $\nabla^{\mathrm{g}}$ be the derivation law of the Levi Civita connection of $g$. The first canonical connection of $(M,\varphi,g)$ is the linear connection having the covariant derivative given by
\[
\nabla^{0}_X Y = \nabla^{\mathrm{g}}_X Y -\frac{1}{2} (\nabla^{\mathrm{g}}_X J_{\varphi}) J_{\varphi} Y, \quad \forall X, Y \in {\mathfrak X} (M).
\]
\end{defin}

The first canonical connection have been studied in a particular case of almost product Riemannian manifolds in \cite{brno}. More precisely, Lemmas 3.10 and 3.12 of that paper prove that $\nabla^{0}$ is adapted to the almost product Riemannian structure and show how to build the set of adapted covariant derivatives starting from $\nabla^{0}$. The proofs detailed there also runs here. These results and Proposition \ref{teor:equivalence-adapted-connections-metric} allow to show the set of natural covariant derivatives on an almost Golden Riemannian manifold by means of the first canonical connection.

\begin{lema}Let $(M,\varphi,g)$ be an almost Golden Riemannian manifold. The first canonical connection $\nabla^{0}$ is adapted to $(\varphi,g)$. Moreover, the set of natural derivation laws of $(\varphi,g)$ is
\[
\{\nabla^{0}+Q\colon Q \in \mathcal L, g(Q(X,Y),Z)+g(Q(X,Z),Y)=0, \forall X, Y, Z \in {\mathfrak X} (M)\}.
\]
\end{lema}

The first canonical connection is not the unique distinguished adapted connection to the almost Golden Riemannian structure which can be introduced with the help of the induced almost product Riemannian structure. There exists another adapted connection previously introduced in \cite{racsam}, attached to an almost product Riemannian structure in the paracomplex case; i.e., $r=s$. However, the most part of the results and proofs showed there are also still valid in the case $r\neq s$.
  
\begin{teor}[{\cite[Theor. 4.4]{racsam}}]
\label{teor:welladapted}
Let $(M,\varphi,g)$ be an almost Golden Riemannian manifold. There exists a unique derivation law $\nabla^{\mathrm{w}}$ on $M$ satisfying $\nabla^{\mathrm{w}} J_{\varphi}=0$, $\nabla^{\mathrm{w}} g =0$ and
\begin{equation}
\label{eq:welladapted}
g(\mathrm{T}^{\mathrm{w}}(X,Y),Z)-g(\mathrm{T}^{\mathrm{w}}(Z,Y),X)= g(\mathrm{T}^{\mathrm{w}}(J_{\varphi}Z,Y), J_{\varphi}X)-g(\mathrm{T}^{\mathrm{w}}(J_{\varphi}X,Y),J_{\varphi}Z), \quad \forall X, Y, Z \in {\mathfrak X} (M),
\end{equation}
where $\mathrm{T}^{\mathrm{w}}$ denotes the torsion tensor of $\nabla^{\mathrm{w}}$. The well adapted connection $\Gamma^{\mathrm{w}}$ of the almost product Riemannian structure $(J_{\varphi},g)$ is the linear connection having the covariant derivative $\nabla^{\mathrm{w}}$.
\end{teor}

The well adapted connection is a linear connection adapted to $(J_{\varphi},g)$ distinguished by the condition \eqref{eq:welladapted} over its torsion tensor. But that condition also can be written by the following equivalent conditions
\[
g(\mathrm{T}^{\mathrm{w}}(J_{\varphi}^+X,Y),J_{\varphi}^+Z)-g(\mathrm{T}^{\mathrm{w}}(J_{\varphi}^+Z,Y),J_{\varphi}^+X)=0, \quad 
g(\mathrm{T}^{\mathrm{w}}(J_{\varphi}^-X,Y),J_{\varphi}^-Z)-g(\mathrm{T}^{\mathrm{w}}(J_{\varphi}^-Z,Y),J_{\varphi}^-X)=0, 
\]
for all $X, Y, Z$ vector fields on $M$ (see \cite[Equations (4.4) and (4.5)]{racsam}). Then the well adapted connection  can be introduced on almost Golden Riemannian manifolds as follows

\begin{teor}
\label{teor:welladapted2}
Let $(M,\varphi,g)$ be an almost Golden Riemannian manifold. There exists a unique derivation law $\nabla^{\mathrm{w}}$ on $M$ adapted to $(\varphi,g)$ satisfying
\begin{eqnarray*}
g(\mathrm{T}^{\mathrm{w}}(X,Y),Z)&=&g(\mathrm{T}^{\mathrm{w}}(Z,Y),X), \quad \forall X, Z \in D_{\phi}, \forall Y\in {\mathfrak X} (M),\\
g(\mathrm{T}^{\mathrm{w}}(X,Y),Z)&=&g(\mathrm{T}^{\mathrm{w}}(Z,Y),X), \quad \forall X, Z \in D_{\bar\phi}, \forall Y\in {\mathfrak X} (M),
\end{eqnarray*}
where $\mathrm{T}^{\mathrm{w}}$ denotes the torsion tensor of $\nabla^{\mathrm{w}}$. The well adapted connection $\Gamma^{\mathrm{w}}$ of the almost Golden Riemannian structure $(\varphi,g)$ is the linear connection having the covariant derivative $\nabla^{\mathrm{w}}$.
\end{teor}

The main difference between Schouten and Vr\u{a}nceanu connections in the case of almost Golden manifolds and the previous ones in the case of almost Golden Riemannian manifolds is that the first canonical and the well adapted connections can be uniquely defined. Besides, the Levi Civita connection and the well adapted connection are particular cases of functorial connections. The theory of functorial connections has been developed to study connections in the very large framework of $G$-structures (see \cite{munoz2} and \cite{valdes}). The difference lies in the Lie group $G$ of each $G$-structure defined over $M$ by almost Golden and almost Golden Riemannian structures, more precisely, in the Lie algebra $\mathfrak g$ of the Lie group $G$. 

In \cite[Theor 8.1]{valdes}, the authors characterize the existence of a functorial connection attached to a $G$-structure by means of a condition over the Lie algebra $\mathfrak g$ of $G$. The next result show sufficient conditions that assure the existence of such a connection.

\begin{teor}[{\cite[Theor.\ 2.4]{munoz2}, \cite[Theor. 2.1]{valdes}}]
\label{teor:sufficient}
 Let $G\subseteq GL(n,\mathbb{R})$ be a Lie group and let $\mathfrak g$ its Lie algebra. If $\mathfrak g^{(1)} = 0$ and $\mathfrak g$ is invariant under matrix transpositions, then there exists the functorial connection given in \cite[Theor 8.1]{valdes}, where $\mathfrak g^{(1)}=\{ S  \in \mathrm{Hom} (\mathbb{R}^{n}, \mathfrak g) \colon  S(v) w - S(w) v = 0,  \forall v, w \in \mathbb{R}^{n} \}$ denotes the first prolongation of the Lie algebra $\mathfrak g$.
\end{teor}

Moreover, the vanishing of the first prolongation of the Lie algebra is an obstruction to its existence (see \cite[Theor. 2.2]{valdes}). The Lie algebras of the orthogonal Lie group $O(n;\mathbb{R})$ and the Lie group $O(r;\mathbb{R})\times O(s;\mathbb{R})$  fullfill the conditions of Theorem \ref{teor:sufficient}, then there exist a functorial connection attached to $O(n;\mathbb{R})$- and $O(r;\mathbb{R})\times O(s;\mathbb{R})$-structures. Such connections are the Levi Civita of a Riemannian manifold (see \cite[Sec. 2]{munoz2}) and the well adapted of an almost product Riemannian manifold (see Theorem \ref{teor:welladapted}). However, the first prolongation of Lie algebra of the Lie group $GL(r;\mathbb{R})\times GL(s;\mathbb{R})$ does not vanish then  $GL(r;\mathbb{R})\times GL(s;\mathbb{R})$-structures; i.e., almost Golden and almost product manifolds do not admit a functorial connection. By this reason Schouten and Vr\u{a}nceanu type connections are the most useful linear connections one can consider on these manifolds.

\section{Integrability of $\varphi$ and integrability of the $G_{(\varphi,g)}$-structure}

The main purpose of \cite{gezer} is to characterize the integrability of almost Golden structures  on almost Golden Riemannian manifolds. Proposition \ref{teor:goldenintegrable} is the most known characterization of the integrability of the polynomial structure $\varphi$ on an almost Golden Riemannian manifold $(M,\varphi,g)$ (see, {\em e.g.}; \cite{gezer}). In the quoted paper, the authors obtain another sufficient condition to the integrability of $\varphi$, \cite[Prop. 2.4]{gezer}, which is not a necessary condition. Moreover, the condition showed in the quoted result is equivalent to $\nabla^{\mathrm{g}}  J_{\varphi}=0$ or $\nabla^{\mathrm{g}} \varphi=0$; i.e., it is equivalent to the fact that the derivation law of the Levi Civita connection is an adapted connection to $(\varphi, g)$ (see \cite[Prop. 2.5]{gezer}).
\medskip 

The first canonical connection of $(M,\varphi,g)$ is an adequate tool to obtain necessary and sufficient conditions to characterize the integrability of the almost Golden structure $\varphi$.

\begin{lema}
\label{teor:toro-nijenhuis-induced}
 Let $(M,\varphi,g)$ be an almost Golden Riemannian manifold. Then we have the following identity
\[
\mathrm{T}^{0}(J_{\varphi}X,J_{\varphi}Y)+\mathrm{T}^{0}(X,Y)=-\frac{1}{2} N_{J_{\varphi}}(X,Y), \quad  \forall X, Y\in {\mathfrak X} (M),
\]
where $\mathrm{T}^{0}$ denotes the torsion tensor of $\nabla^{0}$.
\end{lema}

{\bf Proof.} By straightforward calculations starting from the definitions of the almost product structure $J_{\varphi}$ induced by $\varphi$ and the first canonical connection $\nabla^{0}$. $\blacksquare$
\medskip

\noindent As a direct consequence of the previous lemma one has the following about the integrability of $\varphi$.

\begin{teor}
\label{teor:toro-varphi}
 Let $(M,\varphi,g)$ be an almost Golden Riemannian manifold. The polynomial structure $\varphi$ is integrable if and only if the torsion tensor of the first canonical connection satisfies
\begin{equation}
\label{eq:toro-varphi}
\mathrm{T}^{0}(X,Y)=0, \quad \forall X, Y \in D_{\phi}, \quad \mathrm{T}^{0}(X,Y)=0, \quad \forall X, Y \in D_{\bar\phi}.
\end{equation}
\end{teor}

{\bf Proof.} Taking into account Lemmas \ref{teor:nijenhuis-varphi-induced} and  \ref{teor:toro-nijenhuis-induced} one has that $\varphi$ is integrable if and only if
\begin{equation}
\label{eq:toro-induced}
\mathrm{T}^{0}(J_{\varphi}X,J_{\varphi}Y)=-\mathrm{T}^{0}(X,Y), \quad \forall X, Y \in {\mathfrak X} (M).
\end{equation}

Let $X,Y$ be vector fields on $M$. Then one has $X=J_{\varphi}^+X + J_{\varphi}^-X$, $Y=J_{\varphi}^+Y+J_{\varphi}^-Y$, where $J_{\varphi}^+X, J_{\varphi}^+Y\in T^+M=D_{\phi}$ and $J_{\varphi}^-X, J_{\varphi}^-Y \in T^-M=D_{\bar\phi}$, then
\begin{eqnarray*}
\mathrm{T}^{0}(X,Y)&=&\mathrm{T}^{0}(J_{\varphi}^+X,J_{\varphi}^+Y)+\mathrm{T}^{0}(J_{\varphi}^+X,J_{\varphi}^-Y)+\mathrm{T}^{0}(J_{\varphi}^-X,J_{\varphi}^+Y)+\mathrm{T}^{0}(J_{\varphi}^-X,J_{\varphi}^-Y),\\
\mathrm{T}^{0}(J_{\varphi}X, J_{\varphi}Y)&=&\mathrm{T}^{0}(J_{\varphi}^+X,J_{\varphi}^+Y)-\mathrm{T}^{0}(J_{\varphi}^+X,J_{\varphi}^-Y)-\mathrm{T}^{0}(J_{\varphi}^-X,J_{\varphi}^+Y)+\mathrm{T}^{0}(J_{\varphi}^-X,J_{\varphi}^-Y).
\end{eqnarray*}
If the torsion tensor $\mathrm{T}^{0}$ satisfies (\ref{eq:toro-varphi}) then the previous identities allow to obtain the condition (\ref{eq:toro-induced}). 
If the torsion tensor $\mathrm{T}^{0}$ satisfies (\ref{eq:toro-induced}) then it is easy to obtain the following equalities
\begin{eqnarray*}
\mathrm{T}^{0}(J_{\varphi}X, J_{\varphi}Y)&=& \mathrm{T}^{0}(X,Y)= -\mathrm{T}^{0}(X,Y), \quad \forall X, Y, \in T^+M=D_{\phi}, \\
\mathrm{T}^{0}(J_{\varphi}X, J_{\varphi}Y)&=& \mathrm{T}^{0}(X,Y)= -\mathrm{T}^{0}(X,Y), \quad \forall X, Y, \in T^-M=D_{\bar\phi},
\end{eqnarray*}
thus the tensor $\mathrm{T}^{0}$ also satisfies the conditions (\ref{eq:toro-varphi}). $\blacksquare$
\medskip

The integrability of a $G$-structure which admits a functorial connection can be characterized with the help of the well adapted connection (see \cite[Theor. 2.3]{valdes}). Then one has
\begin{prop}
\label{teor:g-structure-integrable}
Let $(M,\varphi,g)$ be an almost Golden Riemannian manifold. The $G_{(\varphi,g)}$-structure defined by $(\varphi,g)$ over $M$ is integrable if and only if the torsion and curvature tensors of the well adapted connection of $(\varphi,g)$ introduced in Theorem \ref{teor:welladapted2} vanish.
\end{prop}

The above results also can be expressed by the coincidence of the different linear connections.  First we recall some identities on an almost product Riemannian manifold $(M,J,g)$ which are direct consequence of the definition of $\nabla^{\mathrm{g}} J$ and the equivalent properties \eqref{eq:g-compatible} and \eqref{eq:g-almostproduct}.

\begin{lema}
\label{teor:tensorNJ}
Let $(M,J,g)$ be an almost product Riemannian manifold. The tensor $\nabla^{\mathrm{g}} J$ satisfies the following relations:
\begin{eqnarray*}
g((\nabla^{\mathrm{g}}_X J)Y,Z)&=& g ((\nabla^{\mathrm{g}}_X J)Z,Y), \\
g((\nabla^{\mathrm{g}}_X J) JY, Z) &=& - g ((\nabla^{\mathrm{g}}_X J) Y, JZ),\\
g((\nabla^{\mathrm{g}}_X J) JY, Z) &=& - g((\nabla^{\mathrm{g}}_X J) JZ, Y), \quad \forall X, Y, Z \in {\mathfrak X} (M).
\end{eqnarray*}
\end{lema}

Then we have

\begin{lema}
\label{teor:no-nijenhuis}
Let $(M,\varphi,J)$ be an almost Golden Riemannian manifold. The following relation holds:
\begin{equation}
\label{eq:toro-nijenhuis}
g(\mathrm{T}^{0}(X,Y),Z)-g(\mathrm{T}^{0}(Z,Y),X)+g(\mathrm{T}^{0}(J_{\varphi}X,Y),J_{\varphi}Z)-g(\mathrm{T}^{0}(J_{\varphi}Z,Y),J_{\varphi}X)=\frac{1}{2}g( N_{J_{\varphi}}(X,Z),Y),
\end{equation}
for all vector fields $X, Y, Z$ on $M$.
\end{lema}

{\bf Proof.} The definition of the first canonical connection implies that
\[
\mathrm{T}^{0}(X,Y)=-\frac{1}{2} (\nabla^{\mathrm{g}}_X J_{\varphi}) J_{\varphi} Y +\frac{1}{2} (\nabla^{\mathrm{g}}_Y J_{\varphi}) J_{\varphi} X, \quad \forall X, Y \in {\mathfrak X} (M).
\]
Taking into account that $(M,J_{\varphi},g)$ is an almost product Riemannian manifold, then the identities of Lemma \ref{teor:tensorNJ} allow to obtain the following relations for any $X, Y, Z$ be vector fields on $M$
\begin{eqnarray*}
g(\mathrm{T}^{0}(X,Y),Z)&=& -\frac{1}{2} g((\nabla^{\mathrm{g}}_X J_{\varphi}) J_{\varphi} Y,Z) +\frac{1}{2} g((\nabla^{\mathrm{g}}_Y J_{\varphi}) J_{\varphi} X,Z)\\ 
                       &=&\frac{1}{2} g((\nabla^{\mathrm{g}}_X J_{\varphi}) J_{\varphi}Z,Y) +\frac{1}{2} g((\nabla^{\mathrm{g}}_Y J_{\varphi}) J_{\varphi} X,Z),\\
-g(\mathrm{T}^{0}(Z,Y),X)&=& \frac{1}{2} g((\nabla^{\mathrm{g}}_Z J_{\varphi}) J_{\varphi} Y,X) -\frac{1}{2} g((\nabla^{\mathrm{g}}_Y J_{\varphi}) J_{\varphi} Z,X)\\ 
                       &=&-\frac{1}{2} g((\nabla^{\mathrm{g}}_Z J_{\varphi}) J_{\varphi}X,Y) +\frac{1}{2} g((\nabla^{\mathrm{g}}_Y J_{\varphi}) J_{\varphi} X,Z),\\
g(\mathrm{T}^{0}(J_{\varphi}X,Y),J_{\varphi}Z)&=& -\frac{1}{2} g((\nabla^{\mathrm{g}}_{J_{\varphi}X} J_{\varphi}) J_{\varphi}Y, J_{\varphi}Z) +\frac{1}{2} g((\nabla^{\mathrm{g}}_Y J_{\varphi}) X,J_{\varphi} Z)\\ 
                       &=&\frac{1}{2} g((\nabla^{\mathrm{g}}_{J_{\varphi} X} J_{\varphi}) Z,Y) -\frac{1}{2} g((\nabla^{\mathrm{g}}_Y J_{\varphi}) J_{\varphi} X,Z), \\
-g(\mathrm{T}^{0}(J_{\varphi} Z,Y),J_{\varphi} X)&=& \frac{1}{2} g((\nabla^{\mathrm{g}}_{J_{\varphi} Z} J_{\varphi}) J_{\varphi}Y, J_{\varphi}X) -\frac{1}{2} g((\nabla^{\mathrm{g}}_Y J_{\varphi}) Z,J_{\varphi}X)\\ 
                       &=&-\frac{1}{2} g((\nabla^{\mathrm{g}}_{J_{\varphi}Z} J_{\varphi}) X,Y) -\frac{1}{2} g((\nabla^{\mathrm{g}}_Y J_{\varphi}) J_{\varphi} X,Z),
\end{eqnarray*}
then, summing up the above equalities and taking into account formula \eqref{eq:nijenhuis-nabla} we obtain the identity (\ref{eq:toro-nijenhuis}). $\blacksquare$

Observe that the first member of the identity of Lemma \ref{teor:no-nijenhuis} vanishes if and only if the derivation law of the first canonical connection satisfies \eqref{eq:welladapted}; i.e., if and only if the first canonical and the well adapted are the same connection. Then, as direct consequence of the previous results one has

\begin{corol} Let $(M,\varphi,g)$ be an almost Golden Riemannian manifold. Then
\begin{enumerate}
\item[$i)$] The first canonical and the well adapted connections of $(\varphi,g)$ coincide if and only if $\varphi$ is integrable.

\item[$ii)$] The Levi Civita and the well adapted connections of $\varphi$ coincide if and only if the Levi Civita connection is adapted to $(\varphi, g)$. In this case, both connections also coincide with the first canonical connection.

\end{enumerate}
\end{corol}

\section*{Acknowledgements}

The third author gratefully thanks to the support from the Post-Doctoral Fellowship of Harish Chandra Research Institute, Department of Atomic Energy, Government of India.

\end{document}